\documentclass[twoside, 11pt]{article}

\title{\bf Selection principles related to $\alpha_i$-properties}

\author{\Large Ljubi\v sa D.R. Ko\v cinac\footnote{Supported, in part,
by MNZ\v ZS RS.}}

\date{}

\usepackage{amsfonts}
\usepackage{amssymb}

\newtheorem{theorem}{{\bf Theorem}}
\newtheorem{lemma}[theorem]{{\bf Lemma}}

\newtheorem{proposition}[theorem]{{\bf Proposition}}
\newtheorem{definition}[theorem]{{\bf Definition}}

\newcommand{\proof}{{\bf Proof. \ }}
\newcommand{\eproof}{{\square}}

\newcommand{\sone}{{\sf S}_1}
\newcommand{\gone}{{\sf G}_1}

\newcommand{\cpx}{{\sf C}_p(X)}

\newcommand{\fplusx}{(2^X,{\sf F}^+)}

\newcommand{\kxvplus}{({\mathbb K}(X),{\sf V}^+)}
\newcommand{\fxvplus}{({\mathbb F}(X),{\sf V}^+)}

\newcommand{\naturals}{{\mathbb N}}

\begin{document}
\maketitle

\begin{abstract}
We investigate selection principles which are motivated by
Arhangel'ski\v i's $\alpha_i$-properties, $i=1,2,3,4$, and their
relations with classical selection principles. It will be shown
that they are closely related to the selection principle $\sone$
and often are equivalent to it.
\end{abstract}

\begin{flushleft}
\vspace{.3cm} {\sf Keywords}: Selection principles,
$\omega$-cover, $k$-cover, $\gamma$-cover, $\gamma_k$-cover,
boundedness of topological groups, upper Fell topology, upper
Vietoris topology.
\end{flushleft}

\section{Introduction}

In this paper we use the usual topological notation and
terminology \cite{engelking} and consider infinite Hausdorff
spaces.

Let us fix some more notation and terminology regarding selection
principles and families of open covers of a topological space
which are necessary for this exposition. For more information in
connection with selection principles we refer the interested
reader to the survey papers \cite{iransurv}, \cite{leccesurvey},
\cite{tsabansurvey}.

\smallskip
Let ${\mathcal A}$ and ${\mathcal B}$ be collections of sets of an
infinite set $X$.

\smallskip
The symbol $\sone({\mathcal A},{\mathcal B})$ denotes the
selection principle:
\begin{quote}
For each sequence $(A_n:n\in\naturals)$ of elements of ${\mathcal
A}$ there is a sequence $(b_n:n\in\naturals)$ such that $b_n\in
A_n$ for each $n\in\naturals$ and $\{b_n:n\in\naturals\}$ is an
element of ${\mathcal B}$.
\end{quote}

When both $\mathcal A$ and $\mathcal B$ are the collection
$\mathcal O$ of open covers of a space $X$, then $\sone(\mathcal
O,\mathcal O)$ defines the classical \emph{Rothberger covering
property} (see \cite{rothberger}).

\smallskip
There is an infinite game, denoted $\gone ({\mathcal A},{\mathcal
B})$, corresponding to $\sone(\mathcal{A},\mathcal{B})$. Two
players, ONE and TWO, play a round for each natural number $n$. In
the $n$--th round ONE chooses a set $A_n\in \mathcal A$ and TWO
responds by an element $b_n$ from $A_n$. A play
$A_1,b_1;\cdots;A_n,b_n;\cdots$ is won by TWO if $\{b_n:n \in
\naturals\} \in \mathcal B$; otherwise, ONE wins.

It is easy to see that if ONE does not have a winning strategy in
the game $\gone({\mathcal A},{\mathcal B})$, then the
corresponding selection hypothesis $\sone({\mathcal A},{\mathcal
B})$ is true. However, the converse implication is not always
true.

\medskip
We introduce now new selection principles. The motivation for
these definitions is the Arhangel'ski\v i definition of
$\alpha_i$-properties, $i=1,2,3,4$, introduced in \cite{arhalpha}.
$\mathcal A$ and $\mathcal B$ are as above.

\begin{definition}
{\rm The symbol $\alpha_i(\mathcal A, \mathcal B)$, $i=1,2,3,4$,
denotes the following selection hypothesis:
\begin{quote}
For each sequence $(A_n:n\in\naturals)$ of infinite elements of
${\mathcal A}$ there is an element $B\in \mathcal B$ such that:

\hskip4mm $\alpha_1(\mathcal A,\mathcal B)$: for each
$n\in\naturals$ the set $A_n \setminus B$ is finite;

\hskip4mm $\alpha_2(\mathcal A,\mathcal B)$: for each
$n\in\naturals$ the set $A_n \cap B$ is infinite;

\hskip4mm $\alpha_3(\mathcal A,\mathcal B)$: for infinitely many
$n\in\naturals$ the set $A_n \cap B$ is infinite;

\hskip4mm $\alpha_4(\mathcal A,\mathcal B)$: for infinitely many
$n\in\naturals$ the set $A_n \cap B$ is nonempty.
\end{quote}
}
\end{definition}
Evidently, if all members of $\mathcal A$ are infinite, then
\[
\alpha_1(\mathcal A,\mathcal B) \Rightarrow \alpha_2(\mathcal
A,\mathcal B) \Rightarrow \alpha_3(\mathcal A,\mathcal B)
\Rightarrow \alpha_4(\mathcal A,\mathcal B)
\]
and
\[
\sone({\mathcal A},{\mathcal B}) \Rightarrow \alpha_4(\mathcal
A,\mathcal B).
\]
However, if $\mathcal A$ contains finite members, then
$\alpha_1(\mathcal A,\mathcal B)$ does not not imply
$\alpha_2(\mathcal A,\mathcal B)$, while $\alpha_3(\mathcal
A,\mathcal B)$ fails (see \cite{boazalpha}).

If for a space $X$ and a point $x\in X$, $\Sigma_x$ denotes the
family of nontrivial sequences in $X$ that converge to $x$, then
$X$ has the Arhangel'ski\v i $\alpha_i$-property, $i=1,2,3,4$,  if
for each $x\in X$ the property $\alpha_i(\Sigma_x,\Sigma_x)$,
$i=1,2,3,4$, holds.

It is known that the four properties $\alpha_i(\Sigma_x,\Sigma_x)$
are different from each other \cite{arhalpha}, \cite{nogura} and
that the same holds in topological groups \cite{dima},
\cite{noguradimatanaka}. However, it was shown in
\cite{marionalpha} that in function spaces $\cpx$ and in some
hyperspaces \cite{alphahyp} the properties $\alpha_2$, $\alpha_3$
and $\alpha_4$ are equivalent to each other and to the
corresponding ${\sf S}_1$ property. We shall see here that for
some classes $\mathcal A$ and $\mathcal B$ the properties
$\alpha_2(\mathcal A, \mathcal B)$, $\alpha_3(\mathcal A,\mathcal
B)$ and $\alpha_4(\mathcal A,\mathcal B)$ are closely related (and
often equivalent) to $\sone(\mathcal A,\mathcal B)$.

Let $X$ be a topological space, $x\in X$, $A\subset X$. Then we
use the following notation.

\begin{itemize}
\item $\mathcal O$: the collection of open covers of $X$;
\item $\Omega$: the collection of $\omega$-covers of $X$;
\item $\mathcal K$: the collection of $k$-covers of $X$;
\item $\Gamma$: the collection of $\gamma$-covers;
\item $\Gamma_k$: the collection of $\gamma_k$-covers;
\item $\Omega_x$: the set $\{A\subset X\setminus\{x\}:x\in\overline{A}\}$;
\item $\Sigma_x$: the set of all nontrivial sequences in $X$ that
converge to $x$.
\end{itemize}

An open cover $\mathcal U$ of a space $X$ is called an
\emph{$\omega$-cover} (a \emph{$k$-cover}) if every finite
(compact) subset of $X$ is contained in a member of $\mathcal U$
and $X$ is not a member of $\mathcal U$ (i.e. we consider
non-trivial covers).

An open cover $\mathcal{U}$ of $X$ is said to be a
\emph{$\gamma$-cover} (\emph{$\gamma_k$-cover}) if it is infinite,
and for each finite (compact) subset $A$ of $X$ the set $\{U\in
\mathcal U: A\nsubseteq U\}$ is finite.

Observe that each infinite subset of a $\gamma$-cover
($\gamma_k$-cover) is still a $\gamma$-cover ($\gamma_k$-cover).
So, we may suppose that such covers are countable. Each finite
(compact) subset of an infinite (non-compact) space belongs to
infinitely many elements of an $\omega$-cover ($k$-cover) of the
space.

Recall that a space $X$ is said to be \emph{$\omega$-Lindel\"of}
(\emph{$k$-Lindel\"of}) if every $\omega$-cover ($k$-cover) of $X$
contains a countable $\omega$-cover ($k$-cover).

\section{General results}
In this section we discuss covering and closure-type properties
$\alpha_i(\mathcal A,\mathcal B)$, $i=2,3,4$, in topological
spaces and identify some classes $\mathcal A$ and $\mathcal B$ for
which these properties are equivalent to $\sone(\mathcal
A,\mathcal B)$.

We have already mentioned that every space $X$ satisfying the
Rothberger covering property $\sone(\mathcal O,\mathcal O)$
satisfies also $\alpha_4(\mathcal O,\mathcal O)$. The real line
$\mathbb R$ satisfies all the properties $\alpha_i(\mathcal
O,\mathcal O)$, $i=2,3,4$, but $\mathbb R$ does not have the
Rothberger property.

\smallskip
However, we have the following result.

\begin{theorem} \label{alpha234-gamma-set} For an $\omega$-Lindel\"of
space $X$ the following are equivalent:
\begin{itemize}
\item[$(1)$] $X$ satisfies $\alpha_2(\Omega, \Gamma)$;
\item[$(2)$] $X$ satisfies $\alpha_3(\Omega, \Gamma)$;
\item[$(3)$] $X$ satisfies $\alpha_4(\Omega, \Gamma)$;
\item[$(4)$] $X$ satisfies $\sone(\Omega,\Gamma)$.
\end{itemize}
\end{theorem}
$\proof$ $(3) \Rightarrow (4)$: Let $(\mathcal U_n:n\in\naturals)$
be a sequence of $\omega$-covers of $X$. Assume that for each $n
\in\naturals$ we have $\mathcal U_n = \{U_{n,m}:m\in\naturals\}$.
For every $n \in\naturals$ define
\[
\mathcal{V}_n = \{U_{1,m_1}\cap\cdots\cap U_{n,m_n}:n < m_1 < m_2
< \cdots < m_n,  U_{i,m_i}\in\mathcal{U}_i, i\le
n\}\setminus\{\emptyset\}.
\]
Then each $\mathcal V_n$ is an $\omega$-cover of $X$. By (3) and
the fact that each infinite subset of a $\gamma$-cover is also a
$\gamma$-cover, there is an increasing sequence $n_1 < n_2 <
\cdots$ in $\naturals$ and a $\gamma$-cover $\mathcal V =
\{V_{n_i}:i\in\naturals\}$ such that for each $i\in \naturals$,
$V_{n_i}\in \mathcal V_{n_i}$. Let for each $i\in \naturals$,
\[
V_{n_i} = U_{1,m_1}\cap\cdots\cap U_{n_i,m_{n_i}}, \, \,
U_{j,m_j}\in \mathcal U_j, \, j\le n_i.
\]
Put $n_0=0$. For each $i \ge 0$ and each $n$ with $n_i < n \le
n_{i+1}$ let $H_n$ be the $n$-th coordinate in the chosen
representation of $V_{n_{i+1}}$:
\[
H_n= U_{n,m_{n_{i+1}}}.
\]
For each $n\in \naturals$, $H_n\in \mathcal U_n$ and the set $
\mathcal H:= \{H_n:n\in\naturals\}$ is a $\gamma$-cover of $X$
because $\mathcal V$ is a refinement of $\mathcal H$, and $X
\notin \mathcal H$. Therefore, $X$ satisfies $\sone(\Omega, \Gamma)$.\\

$(4) \Rightarrow (1)$: Let $(\mathcal U_n: n\in\naturals)$ be a
sequence of $\omega$-covers of $X$ and let for each $n \in
\naturals$, $\mathcal U_n = \{U_{n,m}:m\in \naturals\}$. We shall
use the fact that $\sone(\Omega, \Gamma)$ is equivalent to ONE has
no winning strategy in the game $\gone(\Omega,\Gamma)$ on $X$
\cite{coc1}. Define the following strategy $\sigma$ for ONE. ONE's
first move is $\sigma (\emptyset) = \mathcal U_1$. Assuming that
the set $U_{1,m_{i_1}} \in \mathcal U_1$ is TWO's response, ONE
plays $\sigma(U_{1,m_{i_1}})$ to be $\mathcal V(1,m_{i_1}) =
\{U_{1,m}:m
> m_{i_1}\}$, still an $\omega$-cover of $X$. If TWO now chooses a set
$U_{1,m_{i_2}} \in \mathcal V(1,m_{i_1})$, ONE plays
$\sigma(U_{1,m_{i_1}},U_{1,m_{i_2}}) = \mathcal V(1,m_{i_2}) =
\{U_{1,m}: m > m_{i_2}\}$ which is still an $\omega$-cover of $X$.
Then TWO chooses a set $U_{1,m_{i_3}} \in
\sigma(U_{1,m_{i_1}},U_{1,m_{i_2}})$. And so on.

[Note: For each $n\in\naturals$ and each $\mathcal U_n$ moves of
ONE form a new sequence of $\omega$-covers and ensure that from
each $\mathcal U_n$ TWO chooses infinitely many elements.]

Since $\sigma$ is not a winning strategy for ONE, consider a
$\sigma$-play
\[
\sigma(\emptyset), U_{1,m_{i_1}}; \sigma(U_{1,m_{i_1}}),
U_{1,m_{i_2}}; \sigma(U_{1,m_{i_1}},U_{1,m_{i_2}}), U_{1,m_{i_3}};
\cdots
\]
lost by ONE. That means that the sequence $\mathcal W$ consisting
of TWO's moves is a $\gamma$-cover of $X$. As it contains
infinitely many elements from each $\mathcal U_n$, $n\in
\naturals$, $\mathcal W$ witnesses for the sequence $(\mathcal U
_n:n\in\naturals)$ that $X$ has property $\alpha_2(\Omega,
\Gamma)$. $\eproof$

\medskip
Similarly to the proof of Theorem \ref{alpha234-gamma-set} we can
prove the following two theorems. For that we use:

\smallskip
(i) $X$ satisfies $\sone(\mathcal K,\Gamma)$ iff ONE has no
winning strategy in the game $\gone(\mathcal K,\Gamma)$ on $X$
(see \cite{kcov1}).

\smallskip
(ii) $X$ satisfies $\sone(\mathcal K,\Gamma_k)$ iff ONE has no
winning strategy in the game $\gone(\mathcal K,\Gamma_k)$ on $X$
(see \cite{gamahyp}).

\begin{theorem} \label{alpha234-k-gamma-set} For a $k$-Lindel\"of
non-compact space $X$, the properties $\alpha_2(\mathcal K,
\Gamma)$, $\alpha_3(\mathcal K, \Gamma)$, $\alpha_4(\mathcal K,
\Gamma)$ and $\sone(\mathcal K,\Gamma)$ are equivalent.
\end{theorem}

\begin{theorem} \label{alpha234-gamma_k-set} For a $k$-Lindel\"of
non-compact space $X$, the properties $\alpha_2(\mathcal K,
\Gamma_k)$, $\alpha_3(\mathcal K, \Gamma_k)$, $\alpha_4(\mathcal
K, \Gamma_k)$ and $\sone(\mathcal K,\Gamma_k)$ are equivalent.
\end{theorem}

\medskip
We also have the following results.

\begin{theorem} \label{alpha234} For a space $X$ and
$\mathcal B\in \{\Gamma, \Gamma_k\}$ the following statements are
equivalent:
\begin{itemize}
\item[$(1)$] $X$ satisfies $\alpha_2(\Gamma_k, \mathcal B)$;
\item[$(2)$] $X$ satisfies $\alpha_3(\Gamma_k, \mathcal B)$;
\item[$(3)$] $X$ satisfies $\alpha_4(\Gamma_k, \mathcal B)$;
\item[$(4)$] $X$ satisfies $\sone(\Gamma_k, \mathcal B)$.
\end{itemize}
\end{theorem}
$\proof$ We have to prove only $(3) \Rightarrow (4)$ and $(4)
\Rightarrow (1)$.

\smallskip
$(3) \Rightarrow (4)$: Let $(\mathcal U_n:n\in\naturals)$ be a
sequence of $\gamma_k$-covers of $X$. Enumerate every $\mathcal
U_n$ bijectively as $\mathcal U_n = \{U_{n,m}:m\in\naturals\}$.
For all $n,m\in\naturals$ define
\[
V_{n,m} = U_{1,m}\cap U_{2,m} \cap \cdots \cap U_{n,m}.
\]
Then for each $n$ the set $\mathcal V_n = \{V_{n,m}:m\in
\naturals\}$ is a $\gamma_k$-cover of $X$, because $\mathcal
U_n$'s are $\gamma_k$-covers. By (4) applied to the sequence
$(\mathcal V_n: n\in\naturals)$ there is an increasing sequence
$n_1 < n_2 < \cdots$ in $\naturals$ and a cover $\mathcal V =
(V_{n_i,m_i}:i\in\naturals) \in \mathcal B$ such that for each
$i\in \naturals$, $V_{n_i,m_i}\in \mathcal V_{n_i}$. Put $n_0=0$.
For each $i \ge 0$, each $j$ with $n_i < j \le n_{i+1}$ and each
$V_{n_{i+1},m_{i+1}} = U_{1,m_{i+1}} \cap \cdots \cap
U_{n_{i+1},m_{i+1}}$ put
\[
H_j= U_{j,m_{i+1}}.
\]
For each $j\in \naturals$, $H_j\in \mathcal U_j$ and the set
$\{H_j:j\in\naturals\}$ is in $\mathcal B$ because this set is
refined by $\mathcal V$ which is in $\mathcal B$. So, $X$
satisfies $\sone(\Gamma_k, \mathcal B)$.\\

$(4) \Rightarrow (1)$: Let $(\mathcal U_n: n\in\naturals)$ be a
sequence of $\gamma_k$-covers of $X$. Suppose that for each
$n\in\naturals$, we have $\mathcal U_n = \{U_{n,m}:m\in
\naturals\}$. Choose an increasing sequence $k_1< k_2 < \cdots <
k_p <\cdots$ of positive integers and for each $n$ and each $k_i$
consider $\mathcal V(n,k_i) := \{U_{n,m}:m \ge k_i\}$. Then each
$\mathcal V(n,k_i)$, $n,i\in\naturals$, is a $\gamma_k$-cover of
$X$. Apply now (1) to the sequence $(\mathcal V(n,k_i):i\in
\naturals, n\in\naturals)$ from $\Gamma_k$ and find a sequence
$(V_{n,k_i}:i,n\in\naturals)$ such that for each $(n,i)\in
\naturals \times \naturals$, $V_{n,k_i}\in \mathcal V(n,k_i)$ and
the set $\mathcal W := \{V_{n,k_i}:n,i \in \naturals\} \in
\mathcal B$. It is easy to see that $\mathcal W$ can be chosen in
such a way that for each $n\in\naturals$ the set $\mathcal U_n
\cap \mathcal W$ is infinite. Therefore, $\mathcal W$ witnesses
for the sequence $(\mathcal U _n:n\in\naturals)$ that $X$ has
property $\alpha_2(\Gamma_k, \mathcal B)$. $\eproof$

\medskip
Notice that in a similar way one can prove that for a space $X$
the properties $\alpha_2(\Gamma, \Gamma)$, $\alpha_3(\Gamma,
\Gamma)$, $\alpha_4(\Gamma, \Gamma)$ and $\sone(\Gamma, \Gamma)$
are equivalent.

\smallskip
As B. Tsaban observed \cite{boazalpha}, the property
$\alpha_1(\Gamma, \Gamma)$ is strictly stronger than $\sone
(\Gamma,\Gamma)$.

\section{Applications to topological groups}

Let $(G,\cdot,\tau)$ be a topological group with the neutral
element $e$ and let $\mathcal B_e$ be a local base at $e$. For
each $U\in \mathcal B_e$ with $U\neq G$ define
\begin{center}
$o(U) = \{x\cdot U:x\in G\}$,\\
$\mathcal O(e) = \{o(U): U\in \mathcal B_e\}$;\\
$\omega(U) = \{F\cdot U:F \in \mathbb F(G)\}$,\\
$ \Omega(e) = \{\omega(U): U\in \mathcal B_e \mbox{ and there is
no $F \in \mathbb F(G)$
with $F\cdot U = G$}\}$;\\
$k(U) = \{K\cdot U:K\in \mathbb K(G)\}$, \\
$\mathcal K(e) = \{k(U): U\in \mathcal B_e \mbox{ and there is no
$K \in \mathbb K(G)$ with $K\cdot U = G$}\}$.
\end{center}
Then clearly, $\mathcal O(e) \subset \mathcal O, \; \Omega(e)
\subset \Omega; \; \mathcal K(e) \subset \mathcal K$.

\medskip
In \cite{coc11} (see also \cite{selunif}, \cite{machuraboaz},
\cite{boaz-o-bound}) Menger-bounded, Rothberger-bounded and
Hurewicz-bounded topological groups have been studied. A
topological group $G$ is \emph{Menger-bounded}
(\emph{Rothberger-bounded}, \emph{Hurewicz-bounded}) if it
satisfies the selection principle $\sone(\Omega(e),\mathcal O)$
($\sone(\mathcal O(e),\mathcal O)$, $\sone(\Omega(e),\Gamma)$).

\medskip
We have the following results. Their proofs are similar, so we
prove only the first of them.

\begin{theorem} \label{alpha4-to-Hur} For a topological group $G$
the following are equivalent:
\begin{itemize}
\item[$(1)$] $G$ satisfies $\alpha_4(\Omega(e), \Gamma)$;
\item[$(2)$] $G$ satisfies $\sone(\Omega(e), \Gamma)$;
\item[$(3)$] $G$ satisfies $\sone(\mathcal K(e), \Gamma)$.
\end{itemize}
\end{theorem}
$\proof$ The implications $(2)\Rightarrow (1)$ and $(2)
\Rightarrow (3)$ are obvious.

\smallskip
$(1) \Rightarrow (2)$: Let $(U_n: n\in\naturals)$ be a sequence of
elements of $\mathcal B_e$. For each $n\in \naturals$ let $V_{n} =
U_{1}\cap U_{2} \cap \cdots \cap U_{n}$ be a member of $\mathcal
B_e$. If we now apply (1) to the sequence $(V_n:n\in\naturals)$ we
find an increasing sequence $n_1 < n_2 < \cdots$ in $\naturals$
and finite sets  $F_{n_i}\subset G$, $i\in\naturals$, so that
$\{F_{n_i}\cdot V_{n_i}:n\in\naturals\}$ is a $\gamma$-cover of
$G$. If $n_0=0$, then for each positive integer $n$ with $n_{i-1}<
n\le n_{i}$, $i\in\naturals$, put $F_n= F_{n_i}$ and $U_n$ to be
the $n$-th component in the representation $U_1\cap \cdots \cap
U_{n_i}$ of $V_{n_i}$. Evidently, $\{F_n\cdot U_n:n\in\naturals\}$
is a $\gamma$-cover of $G$, i.e. the sequence
$(F_n:n\in\naturals)$ guaranties for $(U_n: n\in\naturals)$ that
$G$ satisfies $\sone(\Omega(e),\Gamma)$.\\

$(3) \Rightarrow (2)$: Let $(U_n: n\in\naturals)$ be a sequence of
elements of $\mathcal B_e$. For each $n$ pick a $V_n\in \mathcal
B_e$ so that $V_n^2\subset U_n$. By (3) choose a sequence
$(K_n:n\in\naturals)$ of compact subsets of $G$ such that
$\{K_n\cdot V_n:n\in\naturals\}$ is a $\gamma$-cover of $G$. Next,
for each $n$ pick a finite set $F_n$ in $G$ such that $K_n\subset
F_n\cdot V_n$. Then for each $n$ we have $K_n\cdot V_n\subset
(F_n\cdot V_n)\cdot V_n \subset F_n\cdot U_n$ and one concludes
that $\{F_n\cdot U_n:n\in \naturals\}$ is a $\gamma$-cover of $G$.
$\eproof$

\medskip
Notice that $\alpha_2(\Omega(e), \Gamma)$ and $\alpha_3(\Omega(e),
\Gamma)$ are also equivalent to the properties listed in the
theorem above.

\begin{theorem} \label{alpha4-to-sonegamak} For a topological group $G$
the following are equivalent:
\begin{itemize}
\item[$(1)$] $G$ satisfies $\alpha_4(\Omega(e), \Gamma_k)$;
\item[$(2)$] $G$ satisfies $\sone(\Omega(e), \Gamma_k)$;
\item[$(3)$] $G$ satisfies $\sone(\mathcal K(e), \Gamma_k)$.
\end{itemize}
\end{theorem}




\section{$\alpha_i(\mathcal A,\mathcal B)$ properties in
hyperspaces}

In this section we consider $\alpha_i(\mathcal A,\mathcal B)$
properties, $i=2,3,4$, in hyperspaces. We begin with some
definitions that we need.

For a (Hausdorff) space $X$ by $2^X$ we denote the family of all
closed subsets of $X$. $\mathbb K(X)$ is the collection of all
non-empty compact subsets of $X$, and $\mathbb F(X)$ denotes the
family of all non-empty finite subsets of $X$. If $A$ is a subset
of $X$ and $\mathcal A$ a family of subsets of $X$, then we write
\begin{center}
$A^{+} =\{F\in 2^X:F\subset A\}, \, \, \mathcal A^+ =\{A^+: A\in
\mathcal A\}$.
\end{center}

\noindent Notice that we use the same symbol $F$ to denote a
closed subset of $X$ and the point $F$ in $2^X$; from the context
it will be clear what $F$ is.

The \emph{upper Fell topology} ${\sf F}^+$ on $2^X$ is the
topology whose base is the collection
\[
\{(K^c)^+:K \in \mathbb K(X)\}\cup \{2^X\},
\]
while the \emph{upper Vietoris topology} ${\sf V}^+$ has basic
sets of the form $U^+$, $U$ open in $X$. It is clear that
$\kxvplus$ and $\fxvplus$ are considered as subspaces of
$(2^X,{\sf V}^+)$.

In \cite{alphahyp} it was shown that in $\fplusx$ each of
Arhangel'ski\v{i}'s $\alpha_2$, $\alpha_3$ and $\alpha_4$
properties is equivalent to $\sone(\Sigma_E,\Sigma_E)$, $E\in
2^X$. We discuss here some other properties. For similar
consideration see \cite{selhyp}, \cite{rezpythyp}.

\begin{theorem} \label{alpha4-sone-hyp} If $X$ is a space whose
all open subspaces are $k$-Lindel\"of and $E\in 2\sp X$, then the
following statements are equivalent:
\begin{itemize}
\item[$(1)$] $\fplusx$ satisfies $\alpha_4(\Omega_E,\Sigma_E)$;
\item[$(2)$] $\fplusx$ satisfies $\sone(\Omega_E,\Sigma_E)$.
\end{itemize}
\end{theorem}
$\proof$ We have to prove only (1) implies (2). Let $(\mathcal
A_n:n\in\naturals)$ be a sequence of elements of $\Omega_E$. Since
each open subspace of $X$ is $k$-Lindel\"of, $\fplusx$ has
countable tightness (see \cite{chv}, \cite{kcov2}) and we may
assume that for each $n\in\naturals$, \, $\mathcal A_n$ is
countable, say $\mathcal A_n = \{A_{n,m}:m\in\naturals\}$. For
each $n$ let $\mathcal B_n$ be the collection of all sets of the
form
\[
A_{1,m_1}\cup A_{2,m_2} \cup \cdots A_{n,m_n}, \,
A_{i,m_i}\in\mathcal A_i, \, i\le n.
\]
Then each $\mathcal B_{n}$ belongs to $\Omega_E$. Apply $(1)$ to
the sequence $(\mathcal B_n: n\in\naturals)$ of elements of
$\Omega_E$. There exist an increasing sequence $n_1 < n_2 <
\cdots$ in $\naturals$ and a sequence $\mathcal B :=
(B_{n_i}:i\in\naturals) \in \Sigma_E$ such that for each $i\in
\naturals$, $B_{n_i}\in \mathcal B_{n_i}$. Put $n_0=0$ and define
the sequence $(S_n:n\in\naturals)$ in the following manner:

\begin{quote}
If $i \ge 0$, then for each $n$ with $n_i < n\le n_{i+1}$ define
$S_n$ to be $A_{n,m_n}$ in the chosen representation of
$B_{n_{i+1}}$.
\end{quote}

\noindent Note that for each $n\in\naturals$, $S_n\in\mathcal A_n$
and evidently the sequence $\mathcal S:=(S_n:n\in\naturals)$ is an
element of $\Sigma_E$. So, $\mathcal S$ is a selector for the
original sequence $(\mathcal A_n:n\in\naturals)$ showing that
$\fplusx$ satisfies (2). $\eproof$

\medskip
In what follows we shall need the following two simple lemmas.
Because their proofs are similar we prove only the first of them.

\begin{lemma}\label{k-cover-omega-cover} For a space $X$ and an open
cover $\mathcal W$ of $\kxvplus$ the following holds: $\mathcal W$
is an $\omega$-cover of $\kxvplus$ if and only if $\mathcal
U(\mathcal W):=\{U\subset X: U \mbox{ is open in } X \mbox{ and }
U^+ \subset W \mbox{ for some } W\in\mathcal W_n\}$ is a $k$-cover
of $X$.
\end{lemma}
$\proof$ Let $\mathcal W$ be an $\omega$-cover of $\kxvplus$ and
let $K$ be a compact subset of $X$. Then there exists $W
\in\mathcal W$ such that $K \in W$ and consequently there is an
open set $U\subset X$ with $K \in U^+\subset W$. It is understood,
$U\in \mathcal U(\mathcal U)$. On the other hand, $K\subset U$,
i.e. $\mathcal U(\mathcal W)$ is a $k$-cover of $X$.

Conversely, let $\mathcal U(\mathcal W)$ be a $k$-cover of $X$ and
let $\{K_1,\cdots,K_m\}$ be a finite subset of $\kxvplus$. Then
$K=\bigcup_{i=1}^mK_i$ is a compact subset of $X$ and thus $K$ is
contained in some $U\in \mathcal U(\mathcal W)$; pick
$W\in\mathcal W$ such that $U^+\subset W$. From $K_i\subset U$ for
each $i\le m$, it follows $\{K_1,\cdots,K_m\} \subset U^+\subset
W$ which just means that $\mathcal W$ is an $\omega$-cover of
$\kxvplus$. $\eproof$

\begin{lemma}\label{omega-cover-omega-cover} For a space $X$ and
an open cover $\mathcal W$ of $\fxvplus$ the following holds:
$\mathcal W$ is an $\omega$-cover of $\fxvplus$ if and only if
$\mathcal U(\mathcal W):=\{U\subset X: U \mbox{ is open in } X
\mbox{ and } U^+ \subset W \mbox{ for some } W\in\mathcal W_n\}$
is an $\omega$-cover of $X$.
\end{lemma}

We use now the last two lemmas to prove the next two propositions.

\begin{proposition} \label{K(X):omegaL-X:kL} A space $X$ is
$k$-Lindelo\"of if and only if $\kxvplus$ is $\omega$-Lindel\"of.
\end{proposition}
$\proof$ Let $X$ be a $k$-Lindel\"of space and let $\mathcal W$ be
an $\omega$-cover of $\kxvplus$. By Lemma
\ref{k-cover-omega-cover} (and notation from that lemma),
$\mathcal U(\mathcal W)$ is a $k$-cover of $X$. Choose a countable
family $\{U_i:i\in\naturals\} \subset \mathcal U(\mathcal W)$
which is a $k$-cover of $X$. For each $i\in\naturals$ pick $W_i\in
\mathcal W$ such that $U_i^+\subset W_i$. Again by Lemma
\ref{k-cover-omega-cover} $\{W_i:i\in\naturals\}\subset \mathcal
W$ is an $\omega$-cover of $\kxvplus$.

Let us show the converse. Let $\mathcal U$ be a $k$-cover of $X$.
It is easy to check that $\mathcal U^+$ is an $\omega$-cover of
$\kxvplus$. Choose a countable collection $\{U_i^+:i\in\naturals\}
\subset \mathcal U^+$ which is an $\omega$-cover of $\kxvplus$.
Then $\{U_i:i\in\naturals\} \subset \mathcal U$ is a $k$-cover of
$X$, i.e. $X$ is a $k$-Lindel\"of space. $\eproof$

\medskip
Similarly, by using Lemma \ref{omega-cover-omega-cover}, one
obtains

\begin{proposition} \label{F(X):omegaL-X:omegaL} A space $X$ is
$\omega$-Lindelo\"of if and only if $\fxvplus$ is
$\omega$-Lindel\"of.
\end{proposition}

\begin{theorem} \label{K(X)} For a $k$-Lindel\"of space $X$ the
following are equivalent:
\begin{itemize}
\item[$(1)$] $\kxvplus$ satisfies $\alpha_2(\Omega,\Gamma)$;
\item[$(2)$] $\kxvplus$ satisfies $\alpha_3(\Omega,\Gamma)$;
\item[$(3)$] $\kxvplus$ satisfies $\alpha_4(\Omega,\Gamma)$;
\item[$(4)$] $\kxvplus$ satisfies $\sone(\Omega,\Gamma)$;
\item[$(5)$] $X$ satisfies $\sone (\mathcal K, \Gamma_k)$.
\end{itemize}
\end{theorem}
$\proof$ $(1) \Rightarrow (2) \Rightarrow (3) \Rightarrow (4)$
hold for any space.\\

$(4) \Rightarrow (1)$:  By Proposition \ref{K(X):omegaL-X:kL} the
space $\kxvplus$ is $\omega$-Lindel\"of. It remains to apply
Theorem \ref{alpha234-gamma-set}. \\

$(4) \Rightarrow (5)$: Let $(\mathcal U_n:n\in\naturals)$ be a
sequence of $k$-covers of $X$. Then $(\mathcal U_n^+:n \in
\naturals)$ is a sequence of $\omega$-covers of $(\mathbb
K(X),{\sf V}^+)$. Indeed, fix $n$ and let $\{K_1,\cdots, K_m\}$ be
a finite subset of $\mathbb K(X)$. Then $K= K_1\cup\cdots \cup
K_m$ is a compact subset of $X$ and thus there is $U\in\mathcal U$
with $K\subset U$. This means that for each $i\le m$, $K_i\subset
U$, i.e. $K_i\in U^+$. Therefore $\{K_1,\cdots,K_m\} \subset U^+$
and $\mathcal U_n$ is an $\omega$-cover of $\mathcal K(X)$. By (4)
for each $n$, choose an element $U_n^+$ in $\mathcal U_n^+$ such
that the set $\mathcal U^+ = \{U_n^+:n\in \naturals\}$ is a
$\gamma$-cover of $(\mathbb K(X), {\sf V}^+)$. We prove that
$\{U_n:n\in\naturals\}$ is a $\gamma_k$-cover of $X$. Let $K$ be a
compact subset of $X$. Then there is $n_0\in\naturals$ such that
for each $n\ge n_0$ we have ${K}\in U_n^+$, hence $K\subset U_n$.
It shows that $\{U_n:n\in\naturals\}$ is really a $\gamma_k$-cover
of $X$, i.e. that (5) holds.\\

$(5) \Rightarrow (4)$: Let $(\mathcal W_n:n\in\naturals)$ be a
sequence of $\omega$-covers of $(\mathbb K(X),{\sf V}^+)$. For
each $n$ let
\[
\mathcal U_n = \{U\subset X: U \mbox{ is open in } X \mbox{ and }
U^+ \subset W \mbox{ for some } W\in\mathcal W_n\}.
\]
By Lemma \ref{k-cover-omega-cover} each $\mathcal U_n$ is a
$k$-cover of $X$. By (5) applied to the sequence $(\mathcal
U_n:n\in\naturals)$ one can find a sequence $(U_n:n\in\naturals)$
such that for each $n\in \naturals$, $U_n\in\mathcal U_n$ and the
set $\mathcal U = \{U_n:n\in \naturals)$ is a $\gamma_k$-cover of
$X$. For each $U_n\in \mathcal U$ pick an element $W_n\in\mathcal
W_n$ so that $U_n^+\subset W_n$. We claim that $\{W_n:n
\in\naturals\}$ is a $\gamma$-cover of $\kxvplus$ and so it
witnesses for $(\mathcal W_n:n\in\naturals)$ that (4) is
satisfied. Let $K\in \mathbb K(X)$. Then there is $n_0$ such that
for each $n\ge n_0$, $K\subset U_n$, i.e. $K \in U_n^+ \subset
W_n$. $\eproof$

\medskip
It is not difficult to verify that in a similar way, using
Proposition \ref{F(X):omegaL-X:omegaL} and Theorem
\ref{alpha234-gamma-set}, one obtains the following theorem.

\begin{theorem} \label{F(X)} For an $\omega$-Lindel\"of space $X$
the following are equivalent:
\begin{itemize}
\item[$(1)$] $\fxvplus$ satisfies $\alpha_2(\Omega,\Gamma)$;
\item[$(2)$] $\fxvplus$ satisfies $\alpha_3(\Omega,\Gamma)$;
\item[$(3)$] $\fxvplus$ satisfies $\alpha_4(\Omega,\Gamma)$;
\item[$(4)$] $\fxvplus$ satisfies $\sone(\Omega,\Gamma)$;
\item[$(5)$] $X$ satisfies $\sone (\Omega, \Gamma)$.
\end{itemize}
\end{theorem}


\vspace{.5cm}
\begin{flushleft}
Ljubi\v{s}a D.R. Ko\v{c}inac \\
Faculty of Sciences and Mathematics \\
University of Ni\v{s} \\
Vi\v segradska 33 \\
18000 Ni\v{s}, Serbia \\
{\sf lkocinac@ptt.yu}
\end{flushleft}

\end{document}